\documentclass[journal]{IEEEtran}
\ifCLASSOPTIONcompsoc
  \usepackage[nocompress]{cite}
\else
  \usepackage{cite}
\fi
\usepackage{cite}
\usepackage{amsmath,amssymb,amsfonts}
\usepackage{algorithmic}
\usepackage{graphicx}
\usepackage{textcomp}
\usepackage{enumerate}

\usepackage{enumitem}
\usepackage{booktabs}
\usepackage{bbm}
\usepackage{epstopdf}
\usepackage{subfigure}
\usepackage{bm}

\newtheorem{theorem}{Theorem}
\newtheorem{lemma}[theorem]{Lemma}
\newtheorem{remark}{Remark}

\DeclareSymbolFont{bbold}{U}{bbold}{m}{n}
\DeclareSymbolFontAlphabet{\mathbbold}{bbold}

\newcommand{\qed}{\hfill $\blacksquare$}

\newcommand{\until}[1]{\{1,\dots, #1\}}

\newcommand{\integernonnegative}{\ensuremath{\mathbb{Z}}_{\ge 0}}

\newcommand{\subscr}[2]{#1_{\textup{#2}}}

\newcommand{\setdef}[2]{\{#1 \; | \; #2\}}

\newcommand\oprocendsymbol{\hbox{$\square$}}
\newcommand\oprocend{\relax\ifmmode\else\unskip\hfill\fi\oprocendsymbol}

\DeclareSymbolFont{bbold}{U}{bbold}{m}{n}
\DeclareSymbolFontAlphabet{\mathbbold}{bbold}

\newcommand{\Prob}{\mathbb{P}}

\newcommand{\mcN}{\mathcal{N}}

\begin{document}
\title{Convergence of the Heterogeneous Deffuant-Weisbuch Model: A Complete Proof and Some Extensions}
\author{Ge Chen, \IEEEmembership{Senior Member, IEEE}, Wei Su, \IEEEmembership{Member, IEEE}, Wenjun Mei, \IEEEmembership{Member, IEEE}, and Francesco Bullo, \IEEEmembership{Fellow, IEEE}
\thanks{G. Chen is with the Academy of Mathematics and Systems Science, Chinese Academy of Sciences (e-mail: chenge@amss.ac.cn). W. Su is with the School of Mathematics and Statistics, Beijing Jiaotong University (e-mail: su.wei@bjtu.edu.cn). W. Mei is with the Department of Mechanics and Engineering Science, Peking University (e-mail: mei@pku.edu.cn). F. Bullo is with the Center of Control, Dynamical-Systems and Computation, University of California at Santa Barbara (e-mail: bullo@ucsb.edu). }
\thanks{The corresponding author is Wenjun Mei (e-mail: mei@pku.edu.cn).}}

\maketitle

\begin{abstract}
The Deffuant-Weisbuch (DW) model is a  well-known bounded-confidence opinion dynamics that has attracted wide interest. Although the heterogeneous DW model has been studied by simulations over $20$ years,
its convergence proof is open. Our previous paper \cite{GC-WS-WM-FB:20} solves the problem for the case of uniform weighting factors greater than or equal to $1/2$, but the general case remains unresolved.

This paper considers the DW model with heterogeneous confidence bounds and heterogeneous (unconstrained) weighting factors and shows that, with probability one, the opinion of each agent converges to a fixed vector. In other words, this paper resolves the convergence conjecture for the heterogeneous DW model.  Our analysis also clarifies how the convergence speed may be arbitrarily slow under certain parameter conditions.
\end{abstract}


\section{Introduction}\label{sec:introduction}
People form their own opinions by exchanging information on specific topics
in daily life. Frequent interaction between individuals promotes the spread
and diffusion of opinions in the social network.  Opinion dynamics mainly
studies the process of the formation, evolution, and propagation of public
opinion through the interaction among a group of agents.  Because of the
wide applications on stock market, institutional change, election forecast,
task allocation, and coordinated control of social, economic and
engineering systems, opinion dynamics has risen rapidly in recent years as
an interdisciplinary subject of sociology, physics, control science,
psychology, computer science and so on \cite{AVP-RT:17,AVP-RT:18,FB:22}.

The study of opinion dynamics stems from two-step communication flow in
\cite{EK-PFL:55} and social power theory in \cite{JRPF:56,FH:59}. A
popular model for repeatedly assembling opinions with a fixed confidence
matrix was introduced by~\cite{MHDG:74}. Since then, various theories have
been developed, including social impact theory \cite{BL:81}, social power
theory \cite{NEF-ECJ:90}, dynamic social impact theory \cite{BL:96}, and
influence network theory \cite{NEF:98}.

In recent years, the bounded confidence (BC) model of opinion dynamics has
gained attention. In the BC model, individuals only accept the opinions of
others if the difference between their opinions is within a certain
threshold (confidence bound). The Deffuant-Weisbuch (DW) model, also known
as the Deffuant model, was proposed by Deffuant \emph{et al.}
\cite{GD-DN-FA-GW:00}. In this model, at each time step, two individuals
are randomly selected and update their opinions if the other's opinion
falls within their confidence bound. The Hegselmann-Krause (HK) model,
established by Hegselmann and Krause, updates an individual's opinion by
averaging all opinions within their confidence bound at each time.

The BC model, although simple in appearance, is challenging to analyze due
to its nonlinear, state-dependent inter-agent interactions. The convergence
of the homogeneous DW model, in which all agents have the same confidence
bound, was proven in \cite{JL:05a} and its convergence rate was
established in \cite{JZ-GC:15}. Variations of the DW model have also been
considered, such as the modified DW model with an infinite number of agents
in \cite{GC-FF:11} and the DW model with each agent able to choose
multiple neighbors to exchange opinions with at each time step in
\cite{JZ-YH:13, doi:10.1137/16M1083268}.
Although the convergence of heterogeneous DW model has been studied by simulations over $20$ years \cite{GW-GD-FA-JN:02},
its mathematical proof is still lacking.  This problem was listed as a conjecture (Conjecture 3.3.2) in \cite{JL:07},
and was also claimed to be an open problem in \cite{Kou-Shi:12, doi:10.1137/19M1296628, Beklaryan-Akopov:21} or unresolved problem in \cite{JL:07b}.
Our previous
paper \cite{GC-WS-WM-FB:20} solves the convergence problem for the
heterogeneous DW model when the weighting factor is not less than $1/2$,
but the case with weighting factors less than $1/2$ is still open. The
convergence proof of the heterogeneous HK model is also an open problem
\cite{AM-FB:11f}, except for the special case when each agent's confidence
bound is either $0$ or $1$ \cite{BC-CW:17}.

This paper resolves the open problem of the convergence of
the heterogeneous Deffuant-Weisbuch (DW) model in high dimensions, by
showing that the opinions of each agent converge to a fixed vector with
probability one. This solution completes the understanding of the
convergence of the heterogeneous DW model. Additionally, our analysis
clarifies how the convergence speed may be arbitrarily slow under certain
parameter conditions.

The paper is organized as follows.  Section~\ref{Mod_sec} introduces the
general heterogeneous DW model and our main results.  Proofs are reported
in Section~\ref{pomr}. Section \ref{simulations} gives some simulations,
while Section \ref{Conclusions} concludes the paper.

\section{Model and  main results}\label{Mod_sec}
\renewcommand{\thesection}{\arabic{section}}

The original DW model, first proposed in \cite{GD-DN-FA-GW:00}, assumes that there are a group of agents
 $\mathcal{N}=\until{n}$ with  $n\geq 3$, and each agent has a time-varying one-dimensional opinion. In this paper we generalize the DW model to the case of high-dimension, which assumes each agent $i\in\mathcal{N}$ has a $d$-dimensional opinion $\mathbf{x}_i(t)\in\mathbb{R}^d (d\geq 1)$ at each discrete time $t\in \mathbb{N}=\{0,1,2,\ldots\}$.  Let
$$
X(t):=(\mathbf{x}_1(t),\ldots,\mathbf{x}_n(t))^\top\in\mathbb{R}^{n\times d}
$$
denote the state matrix. For each agent $i$, let $r_{i}>0$ denote its \emph{confidence bound}, and $\mu_i\in(0,1)$ denote its weighting factor. We remark that each agent has the same weighting factor in the original DW model. We let $\mathbbm{1}_{\{\cdot\}}$ denote the indicator function, i.e., we let $\mathbbm{1}_{\{\omega\}}=1$ if the property $\omega$ holds true and $\mathbbm{1}_{\{\omega\}}=0$ otherwise. At each time $t\in\integernonnegative$, a pair $\{i_t,j_t\}$ is independently and uniformly selected from the set of all pairs
$$
\mathcal{A}:=\setdef{\{i,j\}}{i,j\in\until{n},i<{j}}.
$$
Let $\|\cdot\|$ denote the $\ell_2$-norm (Euclidean norm). The opinions of the agents $i_t$ and $j_t$ are updated according to
\begin{equation}\label{m1}
\left\{
\begin{aligned}
  \mathbf{x}_{i_t}&(t+1) = \mathbf{x}_{i_t}(t)\\
  &+\mu_{i_t} \mathbbm{1}_{\{\|\mathbf{x}_{j_t}(t)-\mathbf{x}_{i_t}(t)\|
    \leq r_{i_t}\}}(\mathbf{x}_{j_t}(t)-\mathbf{x}_{i_t}(t)),\\
  \mathbf{x}_{j_t}&(t+1) = \mathbf{x}_{j_t}(t)\\
  &+\mu_{j_t} \mathbbm{1}_{\{\|\mathbf{x}_{j_t}(t)-\mathbf{x}_{i_t}(t)\|
    \leq r_{j_t}\}}(\mathbf{x}_{i_t}(t)-\mathbf{x}_{j_t}(t)),
\end{aligned}\right.
\end{equation}
whereas the other agents' opinions remain unchanged:
\begin{gather}\label{m2}
  \mathbf{x}_{k}(t+1)= \mathbf{x}_{k}(t), \enspace \text{for } k\in\until{n}\setminus\{i_t,j_t\}.
\end{gather}
Let $r_{\min}$ and $r_{\max}$ denote the minimal and maximal confidence bound respectively, while
$\mu_{\min}$ and  $\mu_{\max}$ denote the minimal and maximal weighting factor respectively.
If $r_{\min}=r_{\max}$, the DW model is called homogeneous, otherwise heterogeneous.

Let $\Omega=\mathcal{A}^{\infty}$ be the sample space, $\mathcal{F}$ be
the Borel $\sigma$-algebra of $\Omega$, and $\Prob$ be the probability measure
on $\mathcal{F}$ associated with the independent uniform distribution, i.e., for any nonempty sets $\mathcal{T}\subseteq \mathbb{N}$ and $\mathcal{A}_t\subseteq{A}$, $t\in\mathcal{T}$,
\begin{equation*}
\Prob\left(\bigcap_{t\in\mathcal{T}}\left\{\{i_t,j_t\}\in \mathcal{A}_t \right\}\right)=\prod_{t\in\mathcal{T}}\frac{|\mathcal{A}_t|}{|\mathcal{A}|}=\prod_{t\in\mathcal{T}}\frac{2|\mathcal{A}_t|}{n(n-1)},
\end{equation*}
where $|\mathcal{A}_t|$ denotes the cardinality of the set $\mathcal{A}_t$.

Then the probability space of the DW model is written as
$(\Omega,\mathcal{F},\Prob)$.
It is worth mentioning that $\omega\in\Omega$
 refers to an infinite sequence of agent pairs selected for opinion update.

We are now ready to state our main result.

\begin{theorem}{\textbf{\textup{(Convergence of general heterogeneous DW model)}}}
  \label{Main_result}
  Consider the general heterogeneous DW model (\ref{m1})-(\ref{m2}).  For
  any initial state, the opinion of each agent converges to a fixed random
  vector with probability one, i.e., there exist random vectors
  $\mathbf{x}_i^*\in \mathbb{R}^d$, $i\in\mathcal{N}$, such that
  \begin{equation*}
    \Prob \left(\lim_{t\rightarrow\infty}\mathbf{x}_i(t)=\mathbf{x}_i^*\right)=1, ~~\forall i\in\mathcal{N}.
  \end{equation*}
  Moreover, for all $i\neq j$, either $\mathbf{x}_i^*=\mathbf{x}_j^*$ or
  $\|\mathbf{x}_i^*-\mathbf{x}_j^*\|\geq \max\{r_i,r_j\}$.
\end{theorem}

\begin{figure}[htb]
\centering
\includegraphics[width=2in]{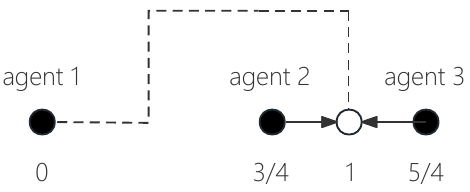}
 \caption{An example where the convergence point is not an equilibrium point in the heterogeneous DW model (\ref{m1})-(\ref{m2}).}\label{remfig1}
 \end{figure}

\begin{remark}
A peculiar fact is that the convergence point
in Theorem \ref{Main_result} may not be an equilibrium point of the system (\ref{m1})-(\ref{m2}). For example, let
$n=3, d=1, (\mathbf{x}_1(0),\mathbf{x}_2(0),\mathbf{x}_3(0))=(0,\frac{3}{4},\frac{5}{4}), (r_1,r_2,r_3)=(\frac{1}{2},\frac{1}{2},1), \mu_1=\mu_2=\mu_3=\frac{1}{3}$.
As shown in Fig. \ref{remfig1}, under the DW protocol it can be calculated that with probability one $\lim_{t\rightarrow\infty}(\mathbf{x}_1(t),\mathbf{x}_2(t),\mathbf{x}_3(t))=(0,1,1)$,
however $(0,1,1)$ is not an equilibrium point because agent $3$ can interact with agent $1$ here.
\end{remark}

Theorem \ref{Main_result} resolves the convergence conjecture for the
heterogeneous DW model and extends the convergence result to the
high-dimensional case with heterogeneous weighting factors.  A natural
question arises regarding the speed of convergence.  Our prior
work~\cite{GC-WS-WM-FB:20} establishes an exponential convergence rate for
the original heterogeneous DW model with the weighting factor being not
less than $1/2$. Here instead we show that this result is not universal and
that the convergence rate can be arbitrarily slow under certain parameter
conditions.

For any $\varepsilon>0$ and any initial state $X(0)$, if $\mathbf{x}_i(t)$ converges
to $\mathbf{x}_i^*$ for any $i\in\mathcal{N}$, we define
\begin{align*}
\tau(\varepsilon)&=\tau_{X(0)}(\varepsilon)\\
&:=\min\left\{t\in\mathbb{Z}^+: \max_{i\in\mathcal{N}}\|\mathbf{x}_i(t)-\mathbf{x}_i^*\|\leq \varepsilon\right\}.
\end{align*}
With this definition we can get the following result.

\begin{theorem}
  \label{result2}
  Consider the general heterogeneous DW model (\ref{m1})-(\ref{m2}). If
  there exist three agents $i,j,k$ satisfying:
  \begin{enumerate}
  \item $r_k\leq r_i< r_j$; and
  \item  $\mu_i<1/2, \mu_j\leq 1/2$, or $\mu_i\leq 1/2,
  \mu_j<1/2$,
  \end{enumerate}
  then with probability one
$$\sup_{X(0)\in\mathbb{R}^{n\times d}} \tau(r_j/2)=\infty.$$
\end{theorem}

The proofs of Theorems~\ref{Main_result}-\ref{result2} are postponed to Section~\ref{pomr}.

\newcommand{\rhomin}{\subscr{\rho}{min}}

\renewcommand{\thesection}{\Roman{section}}
\section{Proofs of Main Results}\label{pomr}
\renewcommand{\thesection}{\arabic{section}}

The proof of Theorem~\ref{Main_result} partially utilizes the approach of our prior work~\cite{GC-WS-WM-FB:20}, which transforms the analysis of a stochastic system into a control algorithm design problem. This method was first introduced to study the original Vicsek model~\cite{GC:17b} and later applied to the noisy HK model~\cite{GC-WS-SD-YH:20} and network evolutionary game \cite{GC-YYY-23}.  Here is the outline the proof of Theorem~\ref{Main_result} using this method:
\begin{itemize}
\item First, we construct a new system called the \emph{DW-control system}
  to aid in the analysis of the DW model. Lemma~\ref{vc} provides the basic
  properties of subsets of agents under the DW-control system.
\item Lemma \ref{merge} demonstrates that two $\varepsilon$-clusters can be
  merged into one $\varepsilon$-cluster in finite time.
\item Finally, Lemmas \ref{robust}-\ref{robust2} prove the key results,
  including the finiteness of the switching times of the DW model topology
  almost surely.
\end{itemize}

It is worth pointing out that the proof of Theorem \ref{Main_result} is  essentially different from our prior work \cite{GC-WS-WM-FB:20}, which heavily depends on
the boundedness of the time for all agents' opinions converging to $\varepsilon$-clusters without connection in DW-control system. This property may not hold when the weighting
factors of some agents is less than $1/2$, see Theorem \ref{result2}
as an example. This paper overcomes this difficulty by utilizing the boundedness of topology switching times in the DW-control system.

\subsection{Some definitions and lemmas}

Consider the DW protocol (\ref{m1})-(\ref{m2}) where, at each time $t$, the pair $\{i_t,j_t\}$ is not selected randomly but instead treated as a control input. In other words, assume that $\{i_t,j_t\}$ is chosen from the set $\mcN$ arbitrarily as a control signal. We call such a control system the \emph{DW-control system}. It is worth mentioning that the DW-control system is an auxiliary tool used only to prove the convergence of the DW model.

For any state $X(t)\in\mathbb{R}^{n\times d}$ at time $t$,
we define the directed edge set $\mathcal{E}_{X(t)}$
by
$$
\mathcal{E}_{X(t)}:=\{(i,j):\|\mathbf{x}_i(t)-\mathbf{x}_j(t)\|\leq r_j\}.
$$
According to this definition, $(i,j)\in\mathcal{E}_{X(t)}$ indicates that agent $i$ can possibly influence agent $j$'s opinion at time $t$.
For any non-empty agent subset $\mathcal{N}^1 \subset \mathcal{N}$, we say  $\mathcal{N}^1$ is a \emph{complete subset} at time $t$
if $\{(i,j),(j,i)\} \subset \mathcal{E}_{X(t)}$  for any $i\neq j\in\mathcal{N}^1$. In other words, complete subsets are vertex subsets that
induce complete undirected subgraphs in the communication graph.
For any two disjoint subsets $\mathcal{N}^1, \mathcal{N}^2\subset \mathcal{N}$, we say
there exist edges between them at time $t$ if there exist agents $i\in\mathcal{N}^1, j\in\mathcal{N}^2$ such that $(i,j)\in \mathcal{E}_{X(t)}$ or $(j,i)\in \mathcal{E}_{X(t)}$;
 otherwise we say there exists no edge between them at time $t$.

For any non-empty agent subset $\mathcal{N}^1 \subseteq \mathcal{N}$, define the diameter of the opinions in $\mathcal{N}^1$ at time $t$ by
\begin{eqnarray*}
\mbox{diam}(\mathcal{N}^1,X(t)):=\max_{i,j\in\mathcal{N}^1}\|\mathbf{x}_i(t)-\mathbf{x}_j(t)\|.
\end{eqnarray*}
If $\lim_{t\rightarrow\infty}\mbox{diam}(\mathcal{N}^1,X(t))=0$,  we say a \emph{consensus} is reached asymptotically among the agents in $\mathcal{N}^1$.

For the DW-control system, we have the following basic properties.

\begin{lemma}\label{vc}
  Consider the DW-control system. Let $T\in\mathbb{Z}^+$ be an integer,
  $\{\{i_{s}',j_{s}'\}\}_{s\geq 0}$ be a sequence of control inputs, and
  $\mathcal{N}^1\subset \mathcal{N}$ be a non-empty subset of agents.
  Then:
  \begin{enumerate}
  \item If there exists no edge between $\mathcal{N}^1$ and
    $\mathcal{N}\setminus\mathcal{N}^1$ at time $t$, then
    $\mbox{diam}(\mathcal{N}^1,X(t+1))\leq
    \mbox{diam}(\mathcal{N}^1,X(t))$.
\item
Assume that during $[t,\infty)$, $\mathcal{N}^1$ is a complete subset and there is no
  edge between $\mathcal{N}^1$ and $\mathcal{N}\setminus\mathcal{N}^1$
.  For $k\in\mathbb{Z}^+$, set $y_k:=1$ if there is one agent communicating its opinion with all
other agents in $\mathcal{N}^1$  during the time $[(k-1)T, kT)$, i.e.,
 there is an
    agent $m_k\in\mathcal{N}^1$ such that
    \begin{equation}\label{vc_0}
      \{l,m_k\}\in \{\{i_{s}',j_{s}'\}\}_{(k-1)T\leq s<kT}, ~~\forall l\in\mathcal{N}^1\backslash\{m_k\},
    \end{equation}
    and $y_k:=0$ otherwise. If $\sum_{k=1}^{\infty} y_k=\infty$,
     then a consensus is reached asymptotically among the agents in $\mathcal{N}^1$.
  \end{enumerate}
\end{lemma}

The proof of Lemma \ref{vc} is postponed to Appendix \ref{proof_vc}.

Let $\varepsilon>0$ be a constant. For any non-empty subset  $\mathcal{N}^1\subseteq \mathcal{N}$, if $\mbox{diam}(\mathcal{N}^1,X(t))\leq \varepsilon$, we say the set $\mathcal{N}^1$ is a \emph{$\varepsilon$-cluster} at time $t$.  It is worth mentioning that a $\varepsilon$-cluster must be a complete subset if $\varepsilon$ is small enough.
Also, we show that two $\varepsilon$-clusters can be merged into one $\varepsilon$-cluster under the DW-control protocol if there are edges between them.

\begin{lemma}\label{merge}
  Let  $X(t)\in\mathbb{R}^{n\times d}, t\geq 0$ be a given state, and
  \begin{eqnarray}\label{dia0_0}
    0<\varepsilon\leq\min_{i\in\mathcal{N}}\min\left\{\frac{r_i}{4}(1-\mu_i)^{ 1+\left\lceil\log_{1-\mu_i}\frac{r_{\min}}{r_i}\right\rceil},\frac{\mu_i r_i}{2}\right\}
  \end{eqnarray}
 be a constant.  Assume that at time $t$ two disjoint $\varepsilon$-clusters
 $\mathcal{N}^1, \mathcal{N}^2 \subset \mathcal{N}$ are connected by one or more edges.
 Then, under the DW-control system, there exists a
 control sequence
 $\{i_{t}',j_{t}'\}$,$\{i_{t+1}',j_{t+1}'\}$,$\ldots$,$\{i_{t+t^*-1}',j_{t+t^*-1}'\}$
 where $t^*$ is bounded by a constant depending on $\varepsilon$ and system
 parameters only, such that $\mathcal{N}^1\cup\mathcal{N}^2$ is a
 $\varepsilon$-cluster at time $t+t^*$, while the opinions of agents
 outside $\mathcal{N}^1\cup\mathcal{N}^2$ keep unchanged.
\end{lemma}

Lemma \ref{merge} seems complex,  however it can be interpreted in an intuitive way:
If two $\varepsilon$-clusters $\mathcal{N}^1$ and $\mathcal{N}^2$ are connected by one or more edges, say between nodes $i\in \mathcal{N}^1$ and
$j\in \mathcal{N}^2$, then we can use a sequence of control inputs to keep bringing $i$ and $j$ close to
each other while alternately ensuring that $\mathcal{N}^1\backslash \{i\}$ remains close to $i$ and $\mathcal{N}^2\backslash \{j\}$ remains
close to $j$, thereby ensuring that $\mathcal{N}^1$ and $\mathcal{N}^2$ eventually merge into a single $\varepsilon$-cluster.
The detailed proof of Lemma \ref{merge} is postponed to Appendix \ref{proof_merge}.

Recall that $\mathbbm{1}_{\{\cdot\}}$ denotes the indicator function.  Let
$$
\xi_{X(t)}:=\sum_{s=t}^{\infty} \mathbbm{1}_{\{\mathcal{E}_{X(s+1)}\neq \mathcal{E}_{X(s)}\}}
$$
denote the number of changes of the edge set $\mathcal{E}_{X(s)}$ during the time interval $[t,\infty)$.
By the repeated use of Lemma \ref{merge}, we can find a sequence of control inputs under the DW-control protocol, such that all agents form some  $\varepsilon$-clusters without edge between them. Based on this result, it can be shown
that $\xi_{X(t)}$ is uniformly upper bounded with positive probability under the DW dynamics.

\begin{lemma}\label{robust}
  Consider the general DW model (\ref{m1})-(\ref{m2}).  Then, there exist
  constants $T\in\mathbb{Z}^+$ and $\delta\in(0,1)$ depending on system
  parameters only, such that for any $t\geq 0$ and
  $X(t)\in\mathbb{R}^{n\times d}$, we have $\Prob\big(\xi_{X(t)}\leq
  T\big)\geq \delta.$
\end{lemma}

Lemma \ref{robust} is only an intermediate result. Next we show that
$\xi_{X(0)}$ is finite with probability one under the DW dynamics.

\begin{lemma}\label{robust2}
For any initial state $X(0)\in\mathbb{R}^{n\times d}$,
 $\Prob\big(\xi_{X(0)}<\infty\big)=1$ under
the general DW model (\ref{m1})-(\ref{m2}).
\end{lemma}

The proofs of Lemmas \ref{robust} and \ref{robust2} are postponed to
Appendices \ref{proof_robust} and \ref{proof_robust2} respectively.

Lemma \ref{robust2} establishes that the underlying topology of the DW model will be fixed after a finite time with probability one.
Combining this with Lemmas \ref{vc} and \ref{merge}, we can further prove that, when the topology becomes fixed, the opinions of agents will form into disjoint
clusters. Finally, the opinions within each cluster will reach consensus, and then Theorem \ref{Main_result} is proved. The detailed proof is put in the next subsection.

\subsection{Proof of Theorem~\ref{Main_result}}

First we show that the dynamic topology of the DW model will almost surely converge to a fixed undirected graph in finite time.
Define $$
\Omega_1:=\{\omega\in\Omega: \xi_{X(0)}(\omega)<\infty\}.
$$
By Lemma \ref{robust2}, we have $\Prob(\Omega_1)=1$.
If $\omega\in\Omega_1$, there exists a finite time $\tilde{\tau}_{X(0)}(\omega)\in [0,\infty)$ such that the edge set remains unchanged when
$s\geq \tilde{\tau}_{X(0)}(\omega)$, i.e.,
\begin{equation}\label{Conr_4}
\mathcal{E}_{X(s+1)}(\omega)=\mathcal{E}_{X(s)}(\omega),~~~~\forall s\geq \tilde{\tau}_{X(0)}(\omega).
\end{equation}

For any  $k\in\mathbb{Z}^+$ and $m\neq l \in\mathcal{N}$ with $r_m\geq r_l$, we set $a_{l,m}^k=a_{m,l}^k:=1$ if the opinion update pair
$\{i_{s},j_{s}\}=\{m,l\}$ for any $$s\in \mathbb{Z}\cap \left[(k-1)\Big\lceil\log_{1-\mu_m} \frac{r_l}{r_m}\Big\rceil,k\Big\lceil\log_{1-\mu_m}\frac{r_l}{r_m}\Big\rceil\right),$$ and $a_{l,m}^k=a_{m,l}^k:=0$ otherwise. Because $\{i_s,j_s\}$ is independently and
uniformly selected from $\mathcal{A}$, we have $\{a_{m,l}^k\}_{k\geq 1}$ is a sequence of independent and identically distributed (i.i.d.) random variables with $\Prob(a_{m,l}^1=1)>0$, and then
\begin{equation}\label{Conr_5}
\Prob\left(\sum_{k=1}^{\infty} a_{m,l}^k=\infty\right)=1,~~\forall m\neq l \in\mathcal{N}.
\end{equation}
Define
$$
\Omega_2:=\left\{\omega\in\Omega: \sum_{k=1}^{\infty} a_{m,l}^k(\omega)=\infty,~~\forall m\neq l \in\mathcal{N}\right\}.
$$
By (\ref{Conr_5}) we get $\Prob(\Omega_2)=1$.

Consider the sample $\omega\in\Omega_1\cap\Omega_2$. For any $(j,i)\in\mathcal{E}_{X(\tilde{\tau}_{X(0)})}(\omega)$ with $r_i\geq r_j$,
by the definition of $\Omega_2$ there exists $k_1\in\mathbb{Z}^+$ such that $a_{i,j}^{k_1}(\omega)=1$ and
$$(k_1-1)\lceil\log_{1-\mu_i} r_j/r_i\rceil\geq \tilde{\tau}_{X(0)}(\omega).$$
By the definition of $a_{i,j}^{k_1}$, $\{i,j\}$ is the opinion update pair during the time $[(k_1-1)\lceil\log_{1-\mu_i} r_j/r_i\rceil,k_1\lceil\log_{1-\mu_i} r_j/r_i\rceil)$,
so with the discussion from (\ref{ml_1}) to (\ref{ml_3_02}) we have $(i,j)\in\mathcal{E}_{X(k_1\lceil\log_{1-\mu_i} r_j/r_i\rceil)}(\omega)$.
Together this with (\ref{Conr_4}), all edges in $\mathcal{E}_{X(s)}(\omega), s\geq\tilde{\tau}_{X(0)}(\omega),$ are undirected, i.e.,
\begin{equation}\label{Conr_6}
(i,j)\in \mathcal{E}_{X(s)}(\omega) \Longleftrightarrow (j,i)\in \mathcal{E}_{X(s)}(\omega).
\end{equation}

Next we show that the time-invariant undirected topology to which the DW model converges, is composed by disjoint complete subgraphs almost surely.
Intuitively, if there is an agent $m$ connecting agents $l$ and $h$ in a time-invariant undirected  topology, then
the opinions of agents $l$ and $h$ must keep getting close to $m$ according to the random interaction rule of the DW model.
 At some point, the opinions of agents  $l$ and $h$ will also get close enough to form an edge. We describe this process
as follows.

Let $K_1$ be a large but fixed integer. For any  $k\in\mathbb{Z}^+$ and $l,m,h \in\mathcal{N}$ with $l<m<h$, we set $b_{l,m,h}^k:=1$ if the opinion update pair $\{i_{s},j_{s}\}=\{l,m\},\{i_{s+1},j_{s+1}\}=\{m,h\},\{i_{s+2},j_{s+2}\}=\{h,l\}$ for all $s\in\{3(k-1)K_1,3(k-1)K_1+3,3(k-1)K_1+6,\ldots,3kK_1-3\}$, and $b_{l,m,h}^k:=0$ otherwise. Define
$$
\Omega_3:=\left\{\omega\in\Omega: \sum_{k=1}^{\infty} b_{l,m,h}^k(\omega)=\infty,~~\forall l<m<h \in\mathcal{N}\right\}.
$$
With the similar discussion to $\Omega_2$ we obtain $\Prob(\Omega_3)=1$.

Consider the sample $\omega\in\Omega_1\cap\Omega_2\cap \Omega_3$. By (\ref{Conr_6}) and (\ref{Conr_4}) all edges in $\mathcal{E}_{X(s)}(\omega)$ are undirected and keep unchanged for any time $s\geq\tilde{\tau}_{X(0)}(\omega)$.

 We assert that for any agents $l'<m'<h'\in\mathcal{N}$, if at time $\tilde{\tau}_{X(0)}(\omega)$ there exist two undirected edges between them, then they must be a complete subset.
 We prove this assertion by contradiction:\\
Assume that the edge set among $\{l',m',h'\}$ is $\{(l',m')$, $(m',l'),(m',h'),(h',m')\}$ for any time $s\geq \tilde{\tau}_{X(0)}(\omega)$ (the other cases can be proved with the similar way).
Since $\omega\in\Omega_3$, there exists an integer $k_2$ satisfying $3(k_2-1)K_1\geq \tilde{\tau}_{X(0)}(\omega)$
and $b_{l',m',h'}^{k_2}(\omega)=1$.
 By the definition of $b_{l',m',h'}^{k_2}$ and the DW protocol (\ref{m1})-(\ref{m2}), we have
\begin{multline}\label{Conr_7}
\left(
  \begin{array}{c}
    \mathbf{x}_{l'}^{\top}(3k_2 K_1) \\
    \mathbf{x}_{m'}^{\top}(3k_2 K_1) \\
    \mathbf{x}_{h'}^{\top}(3k_2 K_1) \\
  \end{array}
\right)=
\left[\left(
  \begin{array}{ccc}
    1 & 0 & 0 \\
    0 & 1-\mu_{m'} & \mu_{m'} \\
    0 & \mu_{h'} & 1-\mu_{h'} \\
  \end{array}
\right)\times\right.\\
\left.\left(
  \begin{array}{ccc}
    1-\mu_{l'} & \mu_{l'} & 0 \\
    \mu_{m'} & 1-\mu_{m'} & 0 \\
    0 & 0 & 1 \\
  \end{array}
\right) \right]^{K_1}\left(
  \begin{array}{c}
    \mathbf{x}_{l'}^{\top}(3(k_2-1) K_1) \\
    \mathbf{x}_{m'}^{\top}(3(k_2-1) K_1) \\
    \mathbf{x}_{h'}^{\top}(3(k_2-1) K_1) \\
  \end{array}
\right)\\
:=Q^{K_1}\left(
  \begin{array}{c}
    \mathbf{x}_{l'}^{\top}(3(k_2-1) K_1) \\
    \mathbf{x}_{m'}^{\top}(3(k_2-1) K_1) \\
    \mathbf{x}_{h'}^{\top}(3(k_2-1) K_1) \\
  \end{array}
\right).
\end{multline}
From straightforward computation it can be checked that $Q$ is a row stochastic, irreducible and aperiodic matrix. Also,
because there are undirected edges between $l'$ and $m'$, and between $m'$ and $h'$  at time $3(k_2-1) K_1$, we have
\begin{multline}\label{Conr_8}
\max_{h_1,h_2\in\{l',m',h'\}} \left\|\mathbf{x}_{h_1}(3(k_2-1)K_1)-\mathbf{x}_{h_2}(3(k_2-1)K_1)\right\|\\
\leq \min\{r_{l'},r_{m'}\}+\min\{r_{h'},r_{m'}\}.
\end{multline}
By (\ref{Conr_7}), (\ref{Conr_8})  and Lemma \ref{app2} we get
\begin{equation*}
\left\|\mathbf{x}_{l'}(3k_2 K_1)-\mathbf{x}_{h'}(3k_2 K_1)\right\|\leq \min\left\{r_{l'},r_{h'}\right\}
\end{equation*}
for large $K_1$, which is contradictory with (\ref{Conr_4}).

From the above assertion, if $\omega\in\Omega_1\cap\Omega_2\cap \Omega_3$ then
the agent set $\mathcal{N}$ can be divided into $K_2=K_2(\omega)$ disjoint complete subsets $\mathcal{N}^1(\omega),\ldots,\mathcal{N}^{K_2}(\omega)$ with no edge between them
 during the time $[\tilde{\tau}_{X(0)}(\omega),\infty)$.

Finally we prove that the agents in each complete subset reach asymptotical consensus almost surely.
Let $K_3\geq n-1$ be an integer. For any $m\in\mathcal{N}$ and $k\in\mathbb{Z}^+$, set $c_m^k:=1$ if
\begin{equation*}
\{l,m\}\in \{\{i_{s}',j_{s}'\}\}_{(k-1)K_3\leq s<kK_3}, ~~\forall l\in\mathcal{N}\backslash\{m\},
\end{equation*}
and $c_m^k:=0$ otherwise. Define
$$\Omega_4:=\left\{\omega\in\Omega: \sum_{k=1}^{\infty} c_{m}^k(\omega)=\infty,~~\forall m \in\mathcal{N}\right\}.$$
With the similar discussion to $\Omega_2$ we can get $\Prob(\Omega_4)=1$.

Consider the case when $\omega\in\cap_{i=1}^4 \Omega_i$. By Lemma \ref{vc} ii), the opinions of all agents in
$\mathcal{N}^j(\omega)$ will converge to a point $\mathbf{z}_j^*(\omega)$ for any $j\in\{1,\ldots,K_2\}$.
Because there is no edge between $\mathcal{N}^i(\omega)$ and $\mathcal{N}^j(\omega)$ with $i\neq j$,
we have  $\|\mathbf{z}_i^*(\omega)-\mathbf{z}_j^*(\omega)\|\geq \max_{k\in \mathbf{z}_i^*(\omega)\cup \mathbf{z}_j^*(\omega)} r_k$.
By the fact of $\Prob(\Omega_i)=1$ for $1\leq i\leq 4$, we have $\Prob(\cap_{i=1}^4 \Omega_i)=1$, then our results are obtained.

\subsection{Proof of Theorem~\ref{result2}}
We first consider the case when the system  has only three agents $i,j,k$.
With the condition $r_i<r_j$, we let $\epsilon>0$ be a constant satisfying
\begin{equation}\label{mr2_1}
\epsilon<(r_j-r_i)\mu_j  ~~\mbox{and}~~ \epsilon\leq \frac{r_i \mu_i \mu_j}{\mu_i+\mu_j}.
\end{equation}
For any large integer $K>0$, set
\begin{eqnarray}\label{mr2_2}
a:=
r_j+\frac{\epsilon[1-(1-\mu_i-\mu_j)^{K}]}{\mu_i+\mu_j(1-\mu_i-\mu_j)^K}<r_j+\frac{\epsilon}{\mu_i},
\end{eqnarray}
where the inequality uses the condition $\mu_i+\mu_j<1$.
Choose
$\mathbf{x}_k(0)=\mathbf{0}$, $\mathbf{x}_i(0)=(r_j-\epsilon/\mu_j,0,\ldots,0)^{\top}$, and
$\mathbf{x}_j(0)=(a,0,\ldots,0)^{\top}$.
By (\ref{mr2_1}) and (\ref{mr2_2}) we have
\begin{eqnarray}\label{m2_3}
\left\{
\begin{aligned}
&3 r_{\max}>a+r_{\max},\\
&a-\left(r_j-\frac{\epsilon}{\mu_j}\right)<r_i,\\
&r_j-\frac{\epsilon}{\mu_j}>r_i.
\end{aligned}\right.
\end{eqnarray}
 Because at the initial time there is edge between agents $i$ and $j$, but no edge between agents $\{i,j\}$ and $k$, so the opinions of agents will change only if the pair $\{i,j\}$ is selected for  opinion update. Let $0\leq t_1<t_2<\cdots$ denote the times when the pair $\{i,j\}$ is selected  for  opinion update. By (\ref{m1})-(\ref{m2}) and the choose of initial state it can be computed that
\begin{eqnarray*}\label{m2_4}
&&\mathbf{x}_j(t_h+1)-\mathbf{x}_i(t_h+1)\\
&&=\left(1-\mu_i-\mu_j\right)^h \left[\mathbf{x}_j(0)-\mathbf{x}_i(0)\right], ~~h=1,2,\ldots,K,
\end{eqnarray*}
and
\begin{eqnarray}\label{m2_5}
&&\mathbf{x}_j(t_{K}+1)\nonumber\\
&&=\mathbf{x}_j(t_{K})-\mu_j \left[\mathbf{x}_j(t_K)-\mathbf{x}_i(t_K)\right]\nonumber\\
&&=\mathbf{x}_j(t_{K-1}+1)-\mu_j \left[\mathbf{x}_j(t_{K-1}+1)-\mathbf{x}_i(t_{K-1}+1)\right]\nonumber\\
&&=\mathbf{x}_j(t_{K-1}+1)-\mu_j \left(1-\mu_i-\mu_j\right)^{K-1} \left[\mathbf{x}_j(0)-\mathbf{x}_i(0)\right]\nonumber\\
&&=\cdots= \mathbf{x}_j(0)-\mu_j \sum_{s=0}^{K-1} \left(1-\mu_i-\mu_j\right)^{s} \left[\mathbf{x}_j(0)-\mathbf{x}_i(0)\right]\nonumber\\
&&=\mathbf{x}_j(0)-\mu_j \frac{1-\left(1-\mu_i-\mu_j\right)^{K}}{\mu_i+\mu_j} \left[\mathbf{x}_j(0)-\mathbf{x}_i(0)\right],
\end{eqnarray}
Substituting the values of $\mathbf{x}_i(0)$ and $\mathbf{x}_j(0)$ into (\ref{m2_5}) we can get
$$
\mathbf{x}_j(t_{K}+1)=(r_j,0,\ldots,0)^{\top},
$$
which indicates there is one edge between agents $k$ and $j$ at time $t_K+1$. By Theorem \ref{Main_result} with probability one the opinions of agents $i,j,k$ converge to a same vector, so with probability one
$$
\tau(r_j/2)\geq t_K+1 \geq K.
$$
Our result can be obtained when we let $K$ grow to infinity.

For the case when $n>3$, we can put the initial opinions of $\mathcal{N}\backslash \{i,j,k\}$ far away from $\{i,j,k\}$, such that the agents in $\mathcal{N}\backslash \{i,j,k\}$ do not affect the
evolution of the opinions of agents  $\{i,j,k\}$, therefore our result still holds.

\begin{figure*}
  \centering
    \includegraphics[width=1\linewidth]{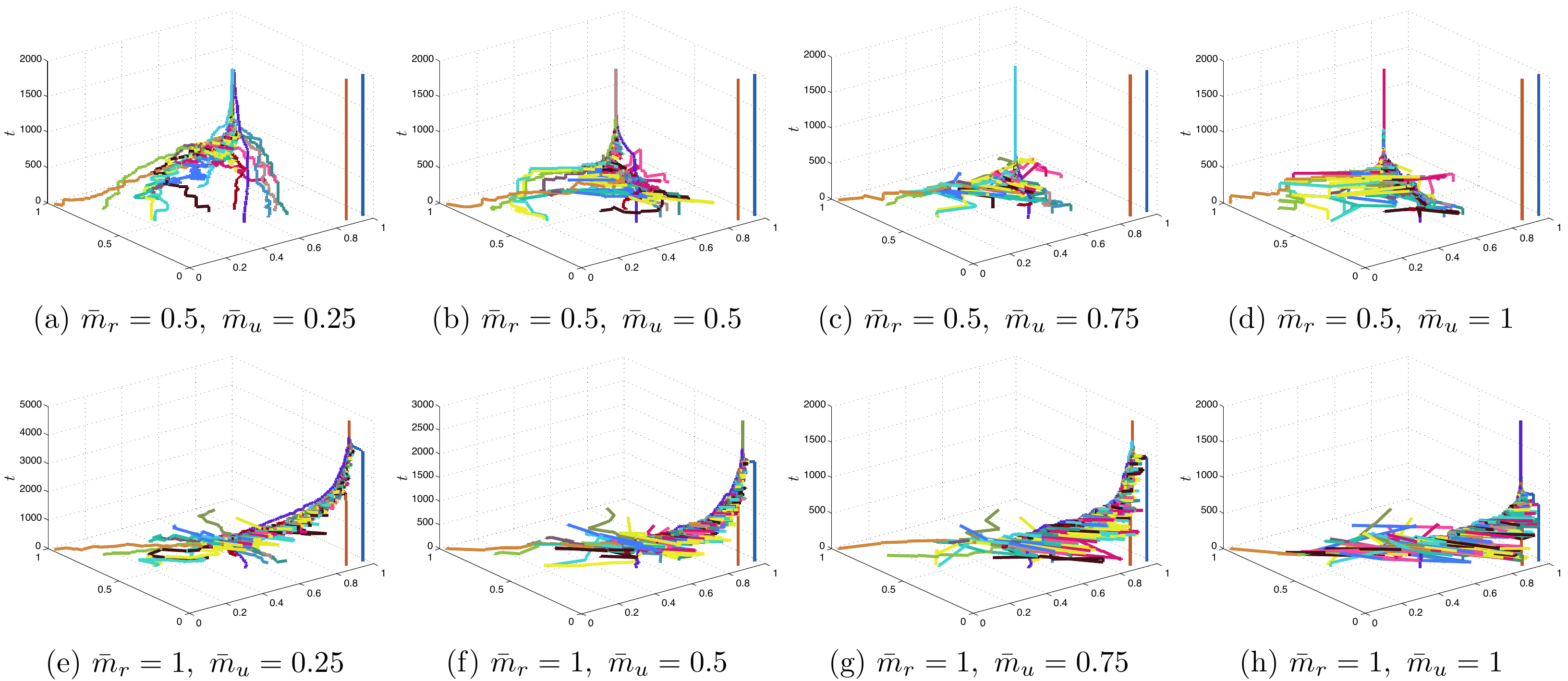}
  \caption{The trajectories of agents' opinions under the system (\ref{m1})-(\ref{m2}) with $n=20$, $d=2$, $\bar{m}_r=0.5,1$, and $\bar{m}_u=0.25,0.5,0.75,1$.}
  \label{Fig1}
\end{figure*}

\begin{figure*}
  \centering
    \includegraphics[width=1\linewidth]{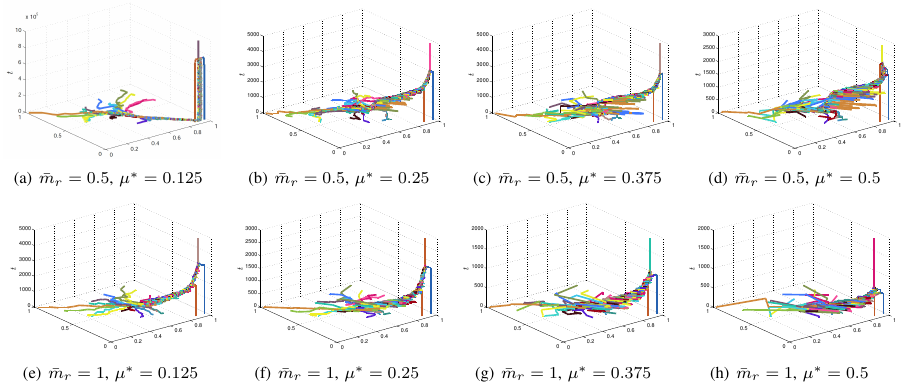}
  \caption{The trajectories of agents' opinions under the system (\ref{m1})-(\ref{m2}) with $n=20$, $d=2$, $\bar{m}_r=0.5,1$, and $\mu^*=0.125,0.25,0.375,5$.}
  \label{Fig2}
\end{figure*}

\renewcommand{\thesection}{\Roman{section}}
\section{Simulations}\label{simulations}
\renewcommand{\thesection}{\arabic{section}}

In this section, we conduct numerical studies on how the dynamic behavior of the heterogeneous DW model depends on its parameters.  We choose $n=20$ and $d=2$, and the initial opinions of all agents are  uniformly and independently distributed in
$[0,1]^2$.
For each agent $i$, its confidence bound $r_i$ and weighting factor $\mu_i$ are set to be $\bar{m}_r \xi_i$ and $\bar{m}_u \zeta_i$
respectively, where $\bar{m}_r$ and $\bar{m}_u$ are two variable parameters, and $\xi_i$ and $\zeta_i$ are two
real numbers randomly selected from $[0,1]$ with independent and uniform distribution. We run the protocol (\ref{m1})-(\ref{m2})
with different $\bar{m}_r$ and $\bar{m}_u$, and the trajectories of the agents' opinions are shown in Fig. \ref{Fig1}.
This figure shows that i) the increase of confidence bounds can promote the consensus of opinions,
however the convergence rate will slow down; ii) the increase of weighting factors can significantly accelerate the convergence rate,
and has no significant impact on the final convergence state of opinions.

We also compare the convergence and convergence speed of the DW model with heterogeneous and homogeneous weighting factors.
Fig. \ref{Fig2} shows the trajectories of the agents' opinions under the protocol (\ref{m1})-(\ref{m2}) with
the weighting factor $\mu_i=\mu^*$ for each agent $i$, while the other configuration is as same as  Fig. \ref{Fig1}.
By comparing Fig.  \ref{Fig1} (a)-(d) and Fig. \ref{Fig2} (a)-(d), it can be seen that with the same expectation of weighting factors ($\mu^*=\bar{m}_u/2$), the  DW model with homogeneous weighting factors
is more likely to reach consensus than with heterogeneous weighting factors. One the other side, by comparing Fig.  \ref{Fig1} (e)-(h) and Fig. \ref{Fig2} (e)-(h),
the DW model with homogeneous and heterogeneous weighting factors has the similar convergence speed when the consensus of opinions is reached.
Moreover, Fig.  \ref{Fig2} (a) indicates that the convergence time may be very large, which verifies Theorem \ref{result2}
in some degree.

\renewcommand{\thesection}{\Roman{section}}
\section{Conclusions}\label{Conclusions}
\renewcommand{\thesection}{\arabic{section}}

The BC opinion dynamics represented by the DW model and HK model has attracted
significant mathematical and sociological attention in recent years. However, due to the mutual coupling between the
underlying topology and agents' states, BC models are generally difficult to analyze. In particular,
the analysis of heterogeneous BC models is very few, even the proof of the most basic convergence property is an open problem that has existed for a long time.
Our previous paper proved the convergence of the heterogeneous DW model with weighting factor not less than $1/2$ \cite{GC-WS-WM-FB:20},
but did not solve the other case. This paper completely proves the convergence of the heterogeneous DW model,
and extends the convergence result to the high-dimensional DW model with heterogenous weighting factors.
Moreover, this paper shows that  the convergence speed may be arbitrarily slow under some parameter conditions.

\appendices
\section{Some Lemmas}

\begin{lemma}\label{app1}
Assume $\mathbf{x},\mathbf{y},\mathbf{z}\in\mathbb{R}^d (d\geq 1)$  satisfying
$\|\mathbf{x}-\mathbf{y}\|=r$, $\|\mathbf{x}-\mathbf{z}\|\leq r$, and $\|\mathbf{y}-\mathbf{z}\|\leq r$.
Let $\mathbf{x}^*=(1-\mu) \mathbf{x}+\mu\mathbf{y}$ with $\mu\in (0,1)$. Then
$\|\mathbf{z}-\mathbf{x}^*\|\leq r\sqrt{1-\mu+\mu^2}$.
\end{lemma}
\textbf{Proof.}
Without loss of generality we assume $\mathbf{x}=\mathbf{0}$, then
\begin{eqnarray}\label{app1_1}
\|\mathbf{z}-\mathbf{x}^*\|^2=\|\mathbf{z}-\mu\mathbf{y}\|^2=\|\mathbf{z}\|^2+\mu^2 r^2 -2\mu \mathbf{y}^{\top}\mathbf{z}.
\end{eqnarray}
By the conditions  $\|\mathbf{y}-\mathbf{z}\|\leq r$ and $\|\mathbf{y}-\mathbf{x}\|=r$ we get
\begin{eqnarray*}\label{app1_2}
r^2\geq \|\mathbf{y}\|^2+\|\mathbf{z}\|^2 -2\mathbf{y}^{\top}\mathbf{z} =r^2+\|\mathbf{z}\|^2 -2\mathbf{y}^{\top}\mathbf{z}\\
~~~~~~~~~~~~~~~ \Longrightarrow  \mu\|\mathbf{z}\|^2 -2\mu\mathbf{y}^{\top}\mathbf{z}\leq 0,
\end{eqnarray*}
then with the condition $\|\mathbf{x}-\mathbf{z}\|\leq r$ we have
\begin{eqnarray*}\label{app1_3}
\mbox{the right-hand side of (\ref{app1_1})}&&\leq (1-\mu)\|\mathbf{z}\|^2+\mu^2 r^2 \\
&&\leq
 (1-\mu)r^2+\mu^2 r^2,
\end{eqnarray*}
which is followed by our result. \qed

\begin{lemma}\label{app2}
Suppose that $P\in[0,1]^{m\times m} (m\geq 2)$ is a row stochastic, irreducible and aperiodic matrix. Let $k, d$ be two positive integers, and  $Y,Z\in\mathbb{R}^{m\times d}$ be two matrices satisfying
$Y=P^k Z$. Then, there exist constants $c>0$ and $\alpha\in(0,1)$ depending on $P$ only such that
\begin{multline*}
\max_{i,j}\|\mbox{Row}_i(Y)-\mbox{Row}_j(Y)\| \\
\leq c \alpha^k \max_{i,j}\|\mbox{Row}_i(Z)-\mbox{Row}_j(Z)\|,
\end{multline*}
where $\mbox{Row}_i(\cdot)$ denotes the $i$-th row of a matrix.
\end{lemma}
\smallskip
\textbf{Proof.}
By Corollary 1.17 in \cite{levin2017markov}, there exists a unique stationary distribution $\bm{\pi}=(\pi_1,\pi_2,\ldots,\pi_m)\in [0,1]^m$ satisfying
$\bm{\pi}=\bm{\pi} P$ and $\sum_{i=1}^m \pi_i=1$. Let $\Pi$ be the matrix with $m$ rows,
each of which is the row vector $\bm{\pi}$. It can be straightforwardly computed that
\begin{equation}\label{napp2_1}
\Pi^2=\Pi,~~P \Pi=\Pi,~~\Pi P=\Pi.
\end{equation}
By the equation $Y=P^k Z$ and repeated use of (\ref{napp2_1}) we have
\begin{eqnarray*}\label{napp2_2}
\begin{aligned}
\Pi Y&=\Pi P^k Z=\Pi P^{k-1} Z\\
&=\cdots= \Pi Z=P\Pi Z=\cdots=P^k \Pi Z,
\end{aligned}
\end{eqnarray*}
and then
\begin{eqnarray}\label{napp2_3}
\begin{aligned}
Y-\Pi Y&=P^k Z- P^k\Pi Z= P^k \left(Z-\Pi Z\right)\\
&= \left(P^k-\Pi\right) \left(Z-\Pi Z\right).
\end{aligned}
\end{eqnarray}
For any $i\in\{1,2,\ldots,m\}$ by (\ref{napp2_3})  we get
\begin{eqnarray}\label{napp2_4}
\begin{aligned}
\mbox{Row}_i(Y)-\bm{\pi} Y&=\mbox{Row}_i(Y-\Pi Y)\\
&= \mbox{Row}_i\left(P^k-\Pi\right) \left(Z-\Pi Z\right)\\
&=\sum_{j=1}^m \left(P^k-\Pi\right)_{ij} \mbox{Row}_j\left(Z-\Pi Z\right).
\end{aligned}
\end{eqnarray}
By Theorem 4.9 in \cite{levin2017markov} there exist constants $c>0$ and $\alpha\in(0,1)$ depending on $P$ only such that
\begin{equation}\label{napp2_5}
\sum_{j=1}^m |\left(P^k-\Pi\right)_{ij}|=\sum_{j=1}^m |\left(P^k\right)_{ij}-\pi_j |\leq \frac{c}{2} \alpha^k.
\end{equation}
On the other hand, with the fact $\sum_{l=1}^m \pi_l=1$ we obtain
\begin{eqnarray}\label{napp2_6}
\begin{aligned}
&\|\mbox{Row}_j\left(Z-\Pi Z\right)\|= \|\mbox{Row}_j(Z)-\bm{\pi}Z\|\\
&= \|\mbox{Row}_j(Z)-[\pi_1 \mbox{Row}_1(Z)+\cdots+\pi_m \mbox{Row}_m(Z)]\|\\
&=\Big\|\sum_{l=1}^m \pi_l \left[\mbox{Row}_j(Z)-\mbox{Row}_l(Z)\right]\Big\|\\
&\leq \max_l \|\mbox{Row}_j(Z)-\mbox{Row}_l(Z)\|,
\end{aligned}
\end{eqnarray}
Substituting (\ref{napp2_5}) and (\ref{napp2_6}) into (\ref{napp2_4}) we have
\begin{equation}\label{napp2_7}
\|\mbox{Row}_i(Y)-\bm{\pi} Y\|\leq \frac{c}{2} \alpha^k \max_{j,l} \|\mbox{Row}_j(Z)-\mbox{Row}_l(Z)\|.
\end{equation}
Since for any $i_1,i_2\in\{1,\ldots,m\}$ we have
\begin{multline*}\label{napp2_8}
\|\mbox{Row}_{i_1}(Y)-\mbox{Row}_{i_2}(Y)\| \\
\leq  \|\mbox{Row}_{i_1}(Y)-\bm{\pi} Y\|+\|\mbox{Row}_{i_2}(Y)-\bm{\pi} Y\|,
\end{multline*}
from (\ref{napp2_7}) our result is obtained.
\qed

\section{Proof of Lemma \ref{vc}}\label{proof_vc}
i) Because there is no edge between $\mathcal{N}^1$ and $\mathcal{N}\setminus\mathcal{N}^1$ at time $t$,
by (\ref{m1})-(\ref{m2}) for each agent $k\in\mathcal{N}^1$, there is a constant $\lambda_k\in [0,1)$ and an agent $k'\in\mathcal{N}^1$ ($k'$ can equal  $k$)
such that
\begin{equation*}\label{vc_1}
  \mathbf{x}_{k}(t+1) = (1-\lambda_k) \mathbf{x}_{k}(t)+\lambda_k\mathbf{x}_{k'}(t).
\end{equation*}
By the property of convex combination we have
\begin{equation*}\label{vc_2}
\mbox{diam}\left(\mathcal{N}^1,X(t+1)\right)\leq \mbox{diam}\left(\mathcal{N}^1,X(t)\right).
\end{equation*}

ii) Without loss of generality we assume $\mathcal{N}^1=\{1,2,\ldots,|\mathcal{N}^1|\}$.
 Define the matrix $P(s)\in[0,1]^{|\mathcal{N}^1|\times |\mathcal{N}^1|}$ as follows:
If $\{i_s',j_s'\}$ is not a subset of $\mathcal{N}^1$, set
$$P(s):=\mbox{diag}(1,\ldots,1)=I_{|\mathcal{N}^1|};$$
 otherwise, set
\begin{eqnarray*}\label{vc_2}
P_{kk}(s):=
\left\{
\begin{aligned}
  1-\mu_{i_s'},~ & \mbox{if } k=i_s'\\
  1-\mu_{j_s'},~ & \mbox{if } k=j_s'\\
  1,~~~~~~~~~~& \mbox{otherwise}
\end{aligned}\right., ~~1\leq k\leq |\mathcal{N}^1|,
\end{eqnarray*}
\begin{eqnarray*}\label{vc_3}
P_{kl}(s):=
\left\{
\begin{aligned}
  \mu_{i_s'},~~ & \mbox{if } (k,l)=(i_s',j_s')\\
  \mu_{j_s'},~~ & \mbox{if } (k,l)=(j_s',i_s')\\
  0,~~~~ & \mbox{otherwise}
\end{aligned}\right., ~~1\leq k\neq j\leq |\mathcal{N}^1|.
\end{eqnarray*}
Thus, $P(s)$ is a row stochastic matrix. Set
\begin{equation}\label{vc_3_1} Z(s):=(\mathbf{x}_1(s),\mathbf{x}_2(s),\ldots,\mathbf{x}_{|\mathcal{N}^1|}(s))^{\top}\in\mathbb{R}^{|\mathcal{N}^1|\times d}.
\end{equation}
Because $\mathcal{N}^1$ is a complete subset, and there is no edge between $\mathcal{N}^1$ and $\mathcal{N}\setminus\mathcal{N}^1$ during $[t,\infty)$,
by the protocol (\ref{m1})-(\ref{m2})
we can get
\begin{equation}\label{vc_4}
Z(s+1)=P(s)Z(s),~~\forall s\geq t.
\end{equation}

For any integer $k\geq \lceil t/T \rceil+1$, let
\begin{equation}\label{vc_6}
\Phi(k):=P(kT-1)P(kT-2)\cdots P((k-1)T),
\end{equation}
then $\Phi(k)$ is still a row stochastic matrix, and by (\ref{vc_4}) we get
\begin{eqnarray}\label{vc_6_1}
\begin{aligned}
Z(kT)&=\Phi(k)Z((k-1)T)=\cdots\\
&=\Phi(k)\Phi(k-1)\cdots\Phi(\lceil t/T \rceil+1)Z(\lceil t/T \rceil T).
\end{aligned}
\end{eqnarray}
Also, by Lemma 2 in \cite{AJ-JL-ASM:02} and the definition of
$P(s)$ we have
\begin{equation}\label{vc_7}
\Phi(k)\geq \delta_T \left[ P(kT-1)+\cdots+P((k-1)T) \right],
\end{equation}
where $\delta_T:=(\frac{(1-\mu_{\max})^2}{2})^{T-1}$ is a constant.

If $y_k=1$, by (\ref{vc_0}), (\ref{vc_7}), and the definition of $P(s)$ we get
\begin{eqnarray}\label{vc_8}
\left\{
\begin{aligned}
&\Phi_{l,m_k}(k)\geq \delta_T \mu_l, ~~\forall l\in\mathcal{N}^1\backslash \{m_k\}\\
&\Phi_{m_k,m_k}(k)\geq \delta_T T (1-\mu_{m_k})
\end{aligned}\right..
\end{eqnarray}
Let $\delta_T':=\min\{\delta_T \mu_{\min},\delta_T T (1-\mu_{\max})\}$, then by (\ref{vc_8}) we have
\begin{multline*}
\min_{h,l}\sum_{m=1}^{|\mathcal{N}^1|} \min\left\{\Phi_{h,m}(k),\Phi_{l,m}(k) \right\}\\
\geq \min_{h,l} \min\left\{\Phi_{h,m_k}(k),\Phi_{l,m_k}(k) \right\}\geq \delta_T'.
\end{multline*}
Thus, by (\ref{vc_6_1}), the fact $\sum_{k=1}^{\infty} y_k=\infty$ and Proposition 5 in \cite{JL:05a},
a consensus is reached asymptotically  among the agents in $\mathcal{N}^1$.

\section{Proof of Lemma \ref{merge}}\label{proof_merge}

\begin{figure}[htb]
\centering
\includegraphics[width=3.3in]{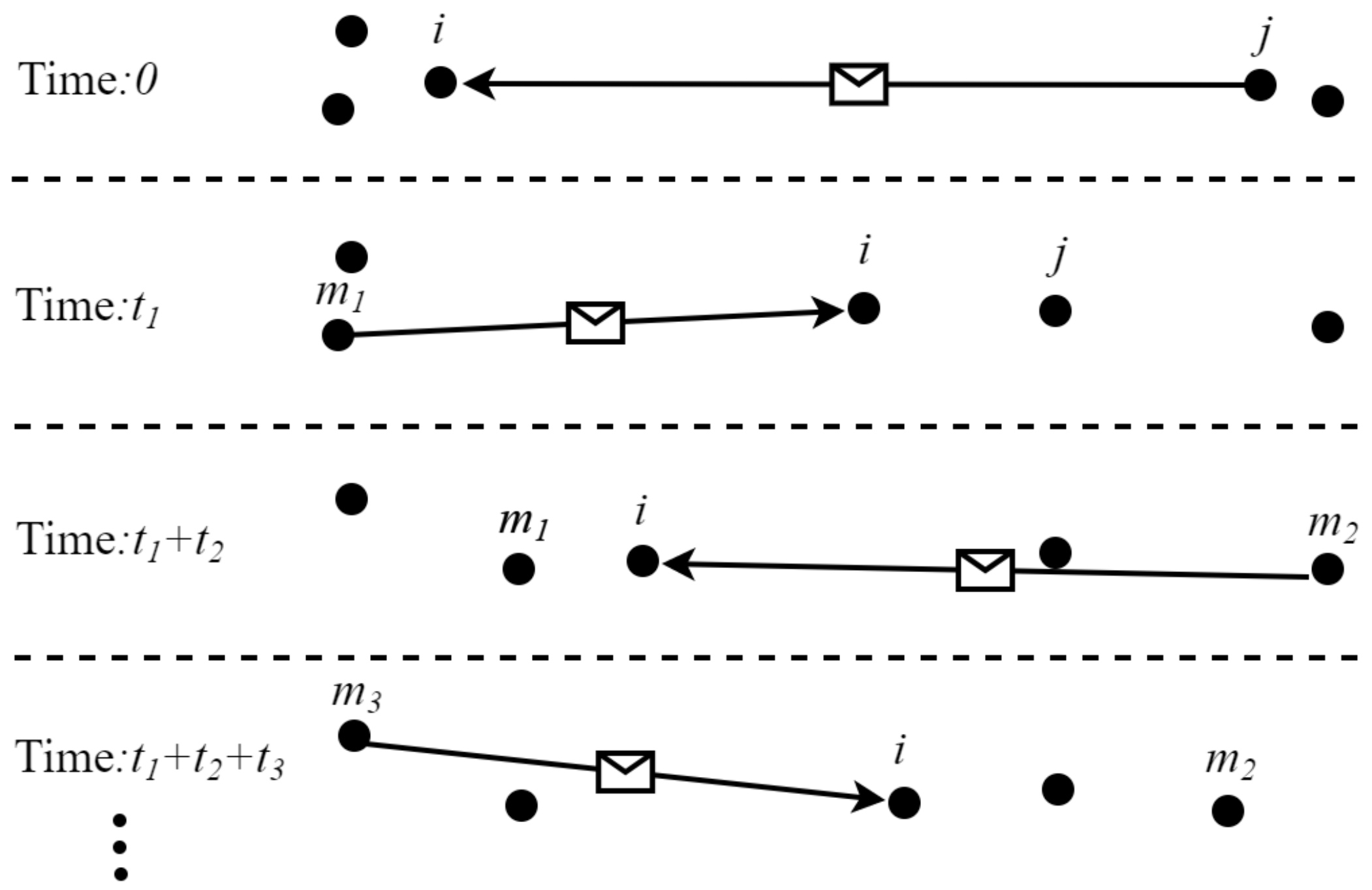}
 \caption{The schematic diagram of the proof of Lemma  \ref{merge}.}\label{Fig4P}
 \end{figure}

Because there are edges between $\mathcal{N}^1$ and $\mathcal{N}^2$ at time $t$, we can find
$i\in\mathcal{N}^1$ and $j\in\mathcal{N}^2$ such that
$$
\|\mathbf{x}_i(t)-\mathbf{x}_j(t)\|\leq \max\{r_i,r_j\}.
$$
Without loss of generality we assume $r_i\geq r_j$. Also, to simplify the exposition we assume $t=0$ and $\mathbf{x}_i(0)=\mathbf{0}$,  then we have $\|\mathbf{x}_j(0)\|\leq r_i$.
We will control agent $i$ to move around and drag other agents together.
 The schematic diagram of this process is shown in Fig. \ref{Fig4P}.

We give the detailed proof by the following six steps.

Step 1: Move agent $i$ towards agent $j$ until agent $j$ also moves once.
Let
\begin{eqnarray}\label{ml_1}
t_1:=
\left\{
\begin{aligned}
 1,~~~~~~~~~~~~~~~~~~~~~~~~~~~~  &\mbox{if}~ \|\mathbf{x}_j(0)\|\leq r_j\\
 1+\left\lceil\log_{1-\mu_i}\frac{r_j}{\|\mathbf{x}_j(0)\|}\right\rceil,~&\mbox{otherwise}
\end{aligned}\right..
\end{eqnarray}
It can be obtained that
\begin{equation}\label{ml_2}
1\leq t_1\leq 1+\left\lceil\log_{1-\mu_i}\frac{r_j}{r_i}\right\rceil\leq 1+\left\lceil\log_{1-\mu_i}\frac{r_{\min}}{r_i}\right\rceil.
\end{equation}

Choose $\{i,j\}$ as the control input at times $0,\ldots,t_1-1$.  Using (\ref{m1}) repeatedly it can be computed that
\begin{eqnarray}\label{ml_3_0}
\left\{
\begin{aligned}
&\mathbf{x}_i(t_1-1)=\left[1-(1-\mu_i)^{t_1-1}\right]\mathbf{x}_j(0),\\
&\mathbf{x}_j(t_1-1)=\mathbf{x}_j(0).
\end{aligned}
\right.
\end{eqnarray}

From this and (\ref{ml_1}) we get
\begin{multline}\label{ml_3_02}
\|\mathbf{x}_i(t_1-1)-\mathbf{x}_j(t_1-1)\|=(1-\mu_i)^{t_1-1} \|\mathbf{x}_j(0)\|\\
\leq (1-\mu_i)^{\log_{1-\mu_i}\frac{r_j}{\|\mathbf{x}_j(0)\|}} \cdot \|\mathbf{x}_j(0)\| =r_j \leq r_i,
\end{multline}
which indicates
\begin{eqnarray}\label{ml_3}
\left\{
\begin{aligned}
&\mathbf{x}_i(t_1)=\left[1-(1-\mu_i)^{t_1}\right]\mathbf{x}_j(0)\\
&\mathbf{x}_j(t_1)=\left[1-\mu_j(1-\mu_i)^{t_1-1}\right]\mathbf{x}_j(0)
\end{aligned}
\right.
\end{eqnarray}
by (\ref{m1}), (\ref{ml_3_0}) and the fact that $\{i,j\}$ is selected to interact opinion at time $t_1-1$.

Step 2:
Reconstruct the coordinates of all individuals' opinions, and show that the ``farthest" agent $m_1$ in $\mathcal{N}^1\cup\mathcal{N}^2$ from agent $i$ is within the confidence bound of agent $i$.
Set
\begin{multline}\label{ml_4}
(\lambda_k(t_1),\delta_k(t_1)):=\\
\left\{
\begin{aligned}
(1-(1-\mu_i)^{t_1} ,\mathbf{0})~~~~~~~~  &\mbox{if}~k=i\\
(1-\mu_j(1-\mu_i)^{t_1-1} ,\mathbf{0})~&\mbox{if}~k=j\\
(0 ,\mathbf{x}_k(0))~~~~~~~~~~~~~~~~~~~~ &\mbox{if}~k\in\mathcal{N}^1 \setminus \{i\}\\
(1 ,\mathbf{x}_k(0)-\mathbf{x}_j(0))~~~~~~~~ &\mbox{if}~k\in\mathcal{N}^2 \setminus \{j\}
\end{aligned}\right.,
\end{multline}
then by (\ref{m2}), (\ref{ml_3}) and the condition diam($\mathcal{N}^i,X(0))\leq \varepsilon$, $i=1,2$,
we have
\begin{eqnarray}\label{ml_5}
\left\{
\begin{aligned}
&\mathbf{x}_k(t_1)=\lambda_k(t_1)\mathbf{x}_j(0)+\delta_k(t_1)\\
&\|\delta_k(t_1)\|\leq \varepsilon
\end{aligned}\right., ~\forall k\in\mathcal{N}^1\cup\mathcal{N}^2.
\end{eqnarray}

Set
\begin{equation}\label{ml_5_1}
\Delta_{\lambda}(t):=\max_{k,l\in\mathcal{N}^1\cup\mathcal{N}^2}|\lambda_{k}(t)-\lambda_l(t)|, ~~\forall t\geq 0.
\end{equation}
Let $m_1$ be the agent (or one of agents) satisfying
\begin{multline}\label{ml_6}
|\lambda_{m_1}(t_1)-\lambda_i(t_1)|\\
=\max_{k\in\mathcal{N}^1\cup\mathcal{N}^2}|\lambda_{k}(t_1)-\lambda_i(t_1)|
\geq \frac{\Delta_{\lambda}(t_1)}{2},
\end{multline}
then
\begin{equation}\label{ml_6_1}
\lambda_{m_1}(t_1)=\max_{k\in\mathcal{N}^1\cup\mathcal{N}^2}\lambda_{k}(t_1) ~\mbox{or}~\min_{k\in\mathcal{N}^1\cup\mathcal{N}^2}\lambda_{k}(t_1).
\end{equation}
By (\ref{ml_4}), (\ref{ml_5}), (\ref{dia0_0}) and (\ref{ml_2})  we have
\begin{eqnarray}\label{ml_7}
\begin{aligned}
&\|\mathbf{x}_{m_1}(t_1)-\mathbf{x}_i(t_1)\|\\
&~\leq |\lambda_{m_1}(t_1)-\lambda_i(t_1)|\times\|\mathbf{x}_j(0)\|+\|\delta_{m_1}(t_1)\|\\
&~\leq \max\left\{(1-\mu_i)^{t_1},1-(1-\mu_i)^{t_1}\right\} r_i + \varepsilon \leq r_i.
\end{aligned}
\end{eqnarray}

Step 3: Move agent $i$ towards  agent $m_1$  until agent $m_1$  also moves once.
Similar to (\ref{ml_1}) we let
\begin{eqnarray*}\label{ml_8}
t_2:=
\left\{
\begin{aligned}
 &1,~~~~~~~~~~~~~~~~~~~~ \mbox{if}~ \|\mathbf{x}_i(t_1)-\mathbf{x}_{m_1}(t_1)\|\leq r_{m_1}&\\
 &1+\left\lceil\log_{1-\mu_i}\frac{r_{m_1}}{ \|\mathbf{x}_i(t_1)-\mathbf{x}_{m_1}(t_1)\|}\right\rceil,~\mbox{otherwise}&
\end{aligned}\right..
\end{eqnarray*}
The inequality (\ref{ml_2}) still holds when we use $t_2$ instead of $t_1$.
Choose $\{i,m_1\}$ as the control input at times $t_1,\ldots,t_1+t_2-1$. With (\ref{ml_7}), similar to (\ref{ml_3}) we can get
\begin{eqnarray}\label{ml_9}
&&\mathbf{x}_i(t_1+t_2)\nonumber\\
&&=\mathbf{x}_{m_1}(t_1)+(1-\mu_i)^{t_2}\left[\mathbf{x}_i(t_1)-\mathbf{x}_{m_1}(t_1) \right]\nonumber\\
&&=\left[\left(1-(1-\mu_i)^{t_2}\right) \lambda_{m_1}(t_1)+(1-\mu_i)^{t_2}\lambda_{i}(t_1) \right] \mathbf{x}_j(0)\nonumber\\
&&~~~+\left(1-(1-\mu_i)^{t_2}\right) \delta_{m_1}(t_1)+(1-\mu_i)^{t_2}\delta_i(t_1)\nonumber\\
&&:=\lambda_i(t_1+t_2) \mathbf{x}_j(0)+\delta_i(t_1+t_2),
\end{eqnarray}
and
\begin{eqnarray}\label{ml_10}
\begin{aligned}
&\mathbf{x}_{m_1}(t_1+t_2)\\
&=\mathbf{x}_{m_1}(t_1)+\mu_{m_1}(1-\mu_i)^{t_2-1}\left[\mathbf{x}_i(t_1)-\mathbf{x}_{m_1}(t_1) \right]\\
&=\left[\left(1-\mu_{m_1}(1-\mu_i)^{t_2-1}\right) \lambda_{m_1}(t_1)\right.\\
&~~~\left.+\mu_{m_1}(1-\mu_i)^{t_2-1}\lambda_{i}(t_1) \right] \mathbf{x}_j(0)\\
&~~~+\left(1-\mu_{m_1}(1-\mu_i)^{t_2-1}\right)\delta_{m_1}(t_1)\\
&~~~+\mu_{m_1}(1-\mu_i)^{t_2-1}\delta_i(t_1)\\
&:=\lambda_{m_1}(t_1+t_2) \mathbf{x}_j(0)+\delta_{m_1}(t_1+t_2).
\end{aligned}
\end{eqnarray}
By (\ref{ml_10}) and (\ref{ml_6}) we get
\begin{eqnarray}\label{ml_10_1}
\begin{aligned}
&|\lambda_{m_1}(t_1+t_2)-\lambda_{m_1}(t_1)|\\
&~=\mu_{m_1}(1-\mu_i)^{t_2-1}\big|\lambda_i(t_1)-\lambda_{m_1}(t_1)\big|\\
&~\geq \mu_{\min}(1-\mu_i)^{\left\lceil\log_{1-\mu_i}\frac{r_{\min}}{r_i}\right\rceil}\frac{\Delta_{\lambda}(t_1)}{2}.
\end{aligned}
\end{eqnarray}

Step 4: Show that all agents in $\mathcal{N}^1\cup\mathcal{N}^2$ are within the confidence bound of agent $i$ at time $t_1+t_2$.
For any  $k\in \mathcal{N}^1\cup\mathcal{N}^2 \setminus \{i,m_1\}$,  let
$\lambda_{k}(t_1+t_2):=\lambda_{k}(t_1)$ and $\delta_{k}(t_1+t_2):=\delta_{k}(t_1)$.
Then, for any  $k\in \mathcal{N}^1\cup\mathcal{N}^2$, by (\ref{ml_5}), (\ref{ml_9}), and (\ref{ml_10})  we have $\|\delta_k(t_1+t_2)\|\leq \varepsilon$,
\begin{eqnarray}\label{ml_10_2}
\mathbf{x}_k(t_1+t_2)=\lambda_k(t_1+t_2)\mathbf{x}_j(0)+\delta_k(t_1+t_2),
\end{eqnarray}
and
\begin{eqnarray}\label{ml_11}
\begin{aligned}
&\|\mathbf{x}_i(t_1+t_2)-\mathbf{x}_{k}(t_1+t_2)\|\\
&~\leq r_i |\lambda_{i}(t_1+t_2)-\lambda_k(t_1+t_2)|+2\varepsilon\\
&~= r_i \big|\lambda_{m_1}(t_1)-\lambda_{k}(t_1+t_2)\\
&~~~+(1-\mu_i)^{t_2}\left(\lambda_i(t_1)-\lambda_{m_1}(t_1)\right)  \big|+2\varepsilon.
\end{aligned}
\end{eqnarray}

If $\lambda_i(t_1)\geq \lambda_{m_1}(t_1)$, by (\ref{ml_6_1}) we have
\begin{eqnarray}\label{ml_12}
\lambda_{m_1}(t_1)=\min_{l\in\mathcal{N}^1\cup\mathcal{N}^2} \lambda_l(t_1).
\end{eqnarray}
Also, because
\begin{eqnarray*}\label{ml_13}
0\leq\min_{l\in\mathcal{N}^1\cup\mathcal{N}^2} \lambda_l(t_1) \leq \lambda_{k}(t_1+t_2)\leq \max_{l\in\mathcal{N}^1\cup\mathcal{N}^2} \lambda_l(t_1)\leq 1,
\end{eqnarray*}
by (\ref{ml_6}), (\ref{ml_12}) and (\ref{dia0_0}) we get
\begin{eqnarray}\label{ml_14}
&&\big|\lambda_{m_1}(t_1)-\lambda_{k}(t_1+t_2)\nonumber\\
&&~~~~~+(1-\mu_i)^{t_2}\left(\lambda_i(t_1)-\lambda_{m_1}(t_1)\right)  \big|\nonumber\\
&&~~\leq \Delta_{\lambda}(t_1)\cdot\max\left\{1-\frac{(1-\mu_i)^{t_2}}{2},(1-\mu_i)^{t_2}  \right\}\nonumber\\
&&~~\leq\max\left\{1-\frac{1}{2}(1-\mu_i)^{1+\left\lceil\log_{1-\mu_i}\frac{r_{\min}}{r_i}\right\rceil},1-\mu_i  \right\}\nonumber\\
&&~~\leq 1-\frac{2\varepsilon}{r_i},
\end{eqnarray}
where $\Delta_{\lambda}(t_1)\leq 1$ can be obtained by  (\ref{ml_4}) and (\ref{ml_5_1}).

If $\lambda_i(t_1)<\lambda_{m_1}(t_1)$, with the similar discussion we can still obtain (\ref{ml_14}).
Substitute (\ref{ml_14}) into (\ref{ml_11}) we have
\begin{eqnarray}\label{ml_19}
\begin{aligned}
&\|\mathbf{x}_i(t_1+t_2)-\mathbf{x}_{k}(t_1+t_2)\|\leq r_i, ~\forall k\in\mathcal{N}^1\cup\mathcal{N}^2,
\end{aligned}
\end{eqnarray}
which indicates that all agents in $\mathcal{N}^1\cup\mathcal{N}^2$ are within the confidence bound of agent $i$ at time $t_1+t_2$.

Step 5: Repeat Steps 3-4 until all agents in $\mathcal{N}^1\cup\mathcal{N}^2$ form a $r_{\min}$-cluster.
Repeat the Steps 3-4 until
\begin{equation}\label{ml_20}
\Delta_{\lambda}(t_1+t_2+\ldots+t_M)\leq (r_{\min}-2\varepsilon)/r_i.
\end{equation}
By (\ref{ml_3}), (\ref{ml_6_1}) and (\ref{ml_10_1}) it can be obtained that
$M$ is bounded by a constant depending on $\mu_k,r_k,k\in\mathcal{N}$, and $\varepsilon$. Let $T:=t_1+t_2+\ldots+t_M$. By (\ref{ml_2}), $T$ is also a finite number.
By (\ref{ml_10_1}), (\ref{ml_20}), and the fact $\|\mathbf{x}_j(0)\|\leq r_i$, we get
\begin{multline}\label{ml_21}
\mbox{diam}(\mathcal{N}^1\cup\mathcal{N}^2,X(T))=\max_{k,l\in \mathcal{N}^1\cup\mathcal{N}^2} \|\mathbf{x}_k(T)-\mathbf{x}_l(T)\|\\
\leq \Delta_{\lambda}(T)\|\mathbf{x}_j(0)\|+2\varepsilon\leq r_{\min},
\end{multline}
which means $\mathcal{N}^1\cup\mathcal{N}^2$ forms a $r_{\min}$-cluster at time $T$.

Step 6: Show that all agents in $\mathcal{N}^1\cup\mathcal{N}^2$ can form a $\varepsilon$-cluster.
For any $s\geq 0$, assume that the agent pair $\{k_s,l_s\}$ has the maximal distance in $\mathcal{N}^1\cup\mathcal{N}^2$ at time $T+s$, i.e.,
\begin{equation}\label{ml_21_1}
\|\mathbf{x}_{k_s}(T+s)-\mathbf{x}_{l_s}(T+s)\|=\mbox{diam}(\mathcal{N}^1\cup\mathcal{N}^2,X(T+s)).
\end{equation}
Choose $\{k_s,l_s\}$ to be control input at time $T+s$, $s\geq 0$, until there exists a time $t^*$ such that $\mbox{diam}(\mathcal{N}^1\cup\mathcal{N}^2,X(t^*))\leq \varepsilon$.
By (\ref{m1}), each element in $\{\mathbf{x}_k(T+s+1)\}_{k\in\mathcal{N}^1\cup\mathcal{N}^2}$ is a convex combination of $\{\mathbf{x}_k(T+s)\}_{k\in\mathcal{N}^1\cup\mathcal{N}^2}$, then
\begin{equation*}\label{ml_22}
\mbox{diam}(\mathcal{N}^1\cup\mathcal{N}^2,X(T+s+1))\leq \mbox{diam}(\mathcal{N}^1\cup\mathcal{N}^2,X(T+s)),
\end{equation*}
and by (\ref{ml_21}),(\ref{ml_21_1}) and (\ref{m1}) we have
\begin{multline}\label{ml_23}
\|\mathbf{x}_{k_s}(T+s+1)-\mathbf{x}_{l_s}(T+s+1)\|\\
=|1-\mu_{k_s}-\mu_{l_s}|\times \mbox{diam}(\mathcal{N}^1\cup\mathcal{N}^2,X(T+s)).
\end{multline}
At the same time, for any agent $m\in\mathcal{N}^1\cup\mathcal{N}^2\setminus\{k_s,l_s\}$, by (\ref{ml_21_1}) and Lemma \ref{app1} we have
\begin{eqnarray}\label{ml_24}
\begin{aligned}
&\left\|\mathbf{x}_{m}(T+s+1)-\mathbf{x}_{k_s}(T+s+1)\right\|\\
&~= \|\mathbf{x}_{m}(T+s)- (1-\mu_{k_s})\mathbf{x}_{k_s}(T+s)\\
&~~~~~~-\mu_{k_s}\mathbf{x}_{l_s}(T+s) \|\\
&~\leq \sqrt{1-\mu_{k_s}+\mu_{k_s}^2}\mbox{diam}\left(\mathcal{N}^1\cup\mathcal{N}^2,X(T+s)\right),
\end{aligned}
\end{eqnarray}
and
\begin{multline}\label{ml_25}
\|\mathbf{x}_{m}(T+s+1)-\mathbf{x}_{l_s}(T+s+1)\|\\
\leq \sqrt{1-\mu_{l_s}+\mu_{l_s}^2}\mbox{diam}\left(\mathcal{N}^1\cup\mathcal{N}^2,X(T+s)\right).
\end{multline}
From (\ref{ml_23}), (\ref{ml_24}) and (\ref{ml_25}) it can be obtained that $t^*$ is bounded by a constant
depending on $\mu_k,r_k,k\in\mathcal{N}$, and $\varepsilon$ only.

By (\ref{m2}),  the opinions of agents outside $\mathcal{N}^1\cup\mathcal{N}^2$ keep unchanged.

\section{Proof of Lemma \ref{robust}}\label{proof_robust}

Without loss of generality we assume $t=0$. Let
  \begin{eqnarray*}
  \varepsilon:=\min_{i\in\mathcal{N}}\min\left\{\frac{r_i}{4}(1-\mu_i)^{ 1+\left\lceil\log_{1-\mu_i}\frac{r_{\min}}{r_i}\right\rceil},\frac{\mu_i r_i}{2}\right\}
  \end{eqnarray*}
be a constant.
Set each agent to be a $\varepsilon$-cluster at the initial time, i.e.,
 $\mathcal{N}_0^1=\{1\},\ldots,\mathcal{N}_0^n=\{n\}$. Let $\mathcal{I}_0:=\{1,2,\ldots,n\}$ be the index set of $\varepsilon$-clusters at the initial time.

If there exists no edge between different $\varepsilon$-clusters at time $0$,
by (\ref{m1})-(\ref{m2}) the opinion of each agent remain unchanged in the time interval $[0,\infty)$,
which indicates $\xi_{X(0)}=0$ certainly. Thus, we only need consider the case when there exist edges
 between different $\varepsilon$-clusters in $\mathcal{I}_0$.

Let $\mathcal{S}_0$ denote the set of index pair $(m,k)$ satisfying i) $m,k\in \mathcal{I}_0,$  and $m<k$;
and ii) there are edges between $\mathcal{N}_0^{m}$ and $\mathcal{N}_0^{k}$ at time $0$.
We record the smallest element in $\mathcal{S}_0$ by the lexicographical order  as $(c_1,c_2)$.
Thus, there exist edges between  $\mathcal{N}_0^{c_1}$ and $\mathcal{N}_0^{c_2}$ at time $0$.
 By Lemma \ref{merge}, under the DW-control protocol, there exists a finite time $t_1$ and
a sequence of control inputs $\{i_0',j_0'\},\{i_{1}',j_{1}'\},\ldots,\{i_{t_1-1}',j_{t_1-1}'\}$
such that $\mathcal{N}_0^{c_1}\cup\mathcal{N}_0^{c_2}$ is a $\varepsilon$-cluster at time $t_1$.
Also, the control inputs of the DW-control protocol is in fact the sample of  the DW system. Thus,
 under the  DW protocol (\ref{m1})-(\ref{m2}), we have
  \begin{equation}\label{robust_1}
    \begin{aligned}
      \Prob&\left(\mbox{diam}(\mathcal{N}_0^{c_1}\cup\mathcal{N}_0^{c_2},X(t_1))\leq \varepsilon\right)\\
      &\qquad\geq \Prob\Big(\bigcap_{s=0}^{t_1-1} \big\{\{i_s,j_s\}=\{i_s',j_s'\}\big\}\Big)\\
      &\qquad=\prod_{s=0}^{t_1-1} \Prob\big(\{i_s,j_s\})=\{i_s',j_s'\}\big)=|\mathcal{A}|^{-t_1}.
    \end{aligned}
  \end{equation}

Let $\mathcal{N}_{t_1}^{c_1}:=\mathcal{N}_{0}^{c_1}\cup \mathcal{N}_{0}^{c_2}$,
$\mathcal{N}_{t_1}^{i}:=\mathcal{N}_{0}^{i}$ for $i\in \mathcal{I}_0 \backslash \{c_1,c_2\}$,
and $\mathcal{I}_{t_1}:=\mathcal{I}_0 \backslash \{c_2\}$.
Define $E_{s}$ to be the event that diam$(\mathcal{N}_{s}^{i},X(s))\leq \varepsilon$
for all $i\in\mathcal{I}_{s}$. By (\ref{robust_1}) and the definition of $\mathcal{N}_{t_1}^{i}$ we have
  \begin{equation}\label{robust_2}
      \Prob\left(E_{t_1}\right)\geq |\mathcal{A}|^{-t_1}.
  \end{equation}

Consider the case when $E_{t_1}$ happens. Define $\tau_1:=\infty$ if there is no edge between different
 $\varepsilon$-clusters in $\mathcal{I}_{t_1}$ during $[t_1,\infty)$, and otherwise
\begin{multline*}
 \tau_1:=\min\{s\geq t_1: \mbox{there are edges between different}\\
 \mbox{ $\varepsilon$-clusters in $\mathcal{I}_{t_1}$ at time $s$}\}.
\end{multline*}
Because there is no edge between different $\varepsilon$-clusters in $\mathcal{I}_{t_1}$ during
$[t_1,\tau_1)$, by Lemma \ref{vc} i) $\mathcal{N}_{\tau_1}^{i}$ is still a $\varepsilon$-cluster for each $i\in \mathcal{I}_{t_1}$,
and also the edge set has no change during $[t_1,\tau_1)$.
Thus, if $\tau_1=\infty$ we can get $\xi_{X(0)}\leq t_1$.

For the case when  $\tau_1<\infty$, similar to $(c_1,c_2)$ we define $(c_3,c_4)$, then
there exists edges between $\mathcal{N}_{t_1}^{c_3}$ and
$\mathcal{N}_{t_1}^{c_4}$ with $c_3,c_4\in\mathcal{I}_{t_1}$ at time $\tau_1$.
By Lemma \ref{merge}, under the DW-control protocol, there exists a finite time $t_2$ and
a sequence of control inputs $\{i_{\tau_1}',j_{\tau_1}'\},\{i_{\tau_1+1}',j_{\tau_1+1}'\},\ldots,\{i_{\tau_1+t_2-1}',j_{\tau_1+t_2-1}'\}$
such that $\mathcal{N}_0^{c_3}\cup\mathcal{N}_0^{c_4}$ forms a $\varepsilon$-cluster at time $\tau_1+t_2$.
Let $\mathcal{N}_{\tau_1+t_2}^{c_3}:=\mathcal{N}_{t_1}^{c_3}\cup \mathcal{N}_{t_1}^{c_4}$,
$\mathcal{N}_{\tau_1+t_2}^{i}:=\mathcal{N}_{t_1}^{i}$ for $i\in \mathcal{I}_{t_1} \backslash \{c_3,c_4\}$,
and $\mathcal{I}_{\tau_1+t_2}:=\mathcal{I}_{t_1} \backslash \{c_4\}$. Similar to (\ref{robust_2}) we can get
  \begin{equation}\label{robust_3}
      \Prob\left(E_{\tau_1+t_2}|E_{t_1},\tau_1<\infty \right)\geq |\mathcal{A}|^{-t_2}.
  \end{equation}

By repeating the above process, we can get a sequence of stop times $\tau_2,\ldots,\tau_{n-2}$ similar to $\tau_1$, and a sequence of finite times $t_3,\ldots,t_{n-1}$ similar to $t_2$,
such that for any $i\in\{2,\ldots,n-2\}$,
  \begin{equation}\label{robust_4}
 \tau_i=\infty  \Longrightarrow \xi_{X(0)}\leq t_1+t_2+\cdots+t_i,
  \end{equation}
  and
    \begin{equation}\label{robust_5}
      \Prob\left(E_{\tau_i+t_{i+1}}|E_{t_1},\tau_1<\infty,\ldots,E_{t_i},\tau_i<\infty \right)
      \geq  |\mathcal{A}|^{-t_{i+1}}.
  \end{equation}
If $E_{\tau_{n-2}+t_{n-1}}$ happens,  by $|\mathcal{I}_{\tau_{n-2}+t_{n-1}}|=1$ we have all agents in $\mathcal{N}$ forms a $\varepsilon$-cluster at time $\tau_{n-2}+t_{n-1}$,
and by Lemma \ref{vc} i) we get
  \begin{equation}\label{robust_6}
\xi_{X(0)}\leq t_1+t_2+\cdots+t_{n-1}.
  \end{equation}
 Let $T:= t_1+t_2+\cdots+t_{n-1}.$ By Lemma \ref{merge}, $T$ can be bounded by a constant depending on system parameters only.
By the discussion from (\ref{robust_2}) to (\ref{robust_6}),  $\xi_{X(0)}\leq T $ if one of the following event happens:
\begin{eqnarray}\label{robust_7}
\begin{aligned}
&F_1:= E_{t_1} \cap \{\tau_1=\infty\},\\
&F_2:= E_{t_1} \cap \{\tau_1<\infty\} \cap  E_{\tau_1+t_2} \cap \{\tau_2=\infty\},\\
&~~~~\cdots\\
&F_{n-2}:= E_{t_1} \cap \{\tau_1<\infty\} \cap \cdots \cap  E_{\tau_{n-3}+t_{n-2}}\\
&~~~~~~~~~~~~~ \cap \{\tau_{n-2}=\infty\},\\
&F_{n-1}:=E_{t_1} \cap \{\tau_1<\infty\} \cap \cdots \cap  E_{\tau_{n-3}+t_{n-2}}\\
&~~~~~~~~~~~~~ \cap \{\tau_{n-2}<\infty\} \cap E_{\tau_{n-2}+t_{n-1}}.
\end{aligned}
\end{eqnarray}
Since $F_1,\ldots,F_{n-1}$ are disjoint events, we have
\begin{eqnarray}\label{robust_8}
&&\Prob\big(\xi_{X(0)}\leq T\big)\geq \Prob\Big(\bigcup_{k=1}^{n-1} F_k \Big)= \sum_{k=1}^{n-1} \Prob\left( F_k \right).
\end{eqnarray}
Let $F_1':= E_{t_1}$, and
\begin{equation*}\label{robust_8_1}
F_{i}':=F_{i-1}' \cap \{\tau_{i-1}<\infty\} \cap E_{\tau_{i-1}+t_{i}}, ~\forall 2\leq i\leq n-1.
\end{equation*}
From (\ref{robust_7}) we have $F_i= F_i' \cap \{\tau_i=\infty\}$ for $1\leq i\leq n-2$, and $F_{n-1}=F_{n-1}'$. By (\ref{robust_5}) we can get
\begin{eqnarray}\label{robust_9}
\Prob\left( F_{n-1}' \right)&&=\Prob\left( F_{n-2}',\tau_{n-2}<\infty, E_{\tau_{n-2}+t_{n-1}}  \right)\nonumber\\
  &&= \Prob\left( F_{n-2}',\tau_{n-2}<\infty \right)\nonumber\\
  &&~~\times\Prob\left(E_{\tau_{n-2}+t_{n-1}}| F_{n-2}',\tau_{n-2}<\infty\right)\nonumber\\
  &&\geq \Prob\left( F_{n-2}',\tau_{n-2}<\infty \right) |\mathcal{A}|^{-t_{n-1}}.
\end{eqnarray}
From (\ref{robust_9}) it can be obtained that
\begin{eqnarray}\label{robust_10}
&&\Prob\left( F_{n-2} \right)+\Prob\left( F_{n-1}' \right)\nonumber\\
&&\geq \Prob\left( F_{n-2}',\tau_{n-2}=\infty \right)+ \Prob\left( F_{n-2}',\tau_{n-2}<\infty \right) |\mathcal{A}|^{-t_{n-1}}\nonumber\\
&&\geq  \Prob\left( F_{n-2}'\right) |\mathcal{A}|^{-t_{n-1}}.
\end{eqnarray}
Similar to (\ref{robust_10}) and with the fact $F_{n-1}=F_{n-1}'$ we have
\begin{eqnarray}\label{robust_11}
&&\Prob\left( F_{n-3} \right)+\Prob\left( F_{n-2} \right)+\Prob\left( F_{n-1} \right)\nonumber\\
&&~\geq \Prob\left( F_{n-3}\right)+ \Prob\left( F_{n-2}' \right) |\mathcal{A}|^{-t_{n-1}}\nonumber\\
&&~\geq \big[\Prob\left( F_{n-3} \right)+ \Prob\left( F_{n-2}' \right)  \big] |\mathcal{A}|^{-t_{n-1}}\nonumber\\
&&~\geq  \Prob\left( F_{n-3}' \right) |\mathcal{A}|^{-t_{n-2}-t_{n-1}},\nonumber\\
&&~\cdots\nonumber\\
&&\Prob\left( F_{1} \right)+\Prob\left( F_{2} \right)+\cdots+\Prob\left( F_{n-1} \right)\nonumber\\
&&~\geq  \Prob\left( F_{1}' \right) |\mathcal{A}|^{-t_2-\cdots-t_{n-1}}\nonumber\\
&&~\geq  |\mathcal{A}|^{-t_1-\cdots-t_{n-1}}=\frac{2^{T}}{n^{T}(n-1)^{T}},
\end{eqnarray}
where the last inequality uses the fact $F_{1}'=E_{t_1}$ and (\ref{robust_2}). Substituting (\ref{robust_11}) into (\ref{robust_8}) yields our result.

\section{Proof of Lemma \ref{robust2}}\label{proof_robust2}
Let $X(0)\in\mathbb{R}^{n\times d}$ be arbitrarily given.
Let $T>0$ be the same constant appearing in Lemma \ref{robust}. Set $T':=T+1$. For any $k\in\mathbb{Z}^+$,
define
\begin{equation*}
\tau_{T'}^k:=\min\Big\{t\in \mathbb{Z}^+: \sum_{s=1}^t \mathbbm{1}_{\{\mathcal{E}_{X(s)}\neq \mathcal{E}_{X(s-1)}\}}= kT'+1\Big\}
\end{equation*}
when $\xi_{X(0)}\geq kT'+1$, and $\tau_{T'}^k:=\infty$ when $\xi_{X(0)}\leq kT'$.
Recall that $\xi_{X(t)}=\sum_{s=t}^{\infty} \mathbbm{1}_{\{\mathcal{E}_{X(s+1)}\neq \mathcal{E}_{X(s)}\}}$. By the definition of $\tau_{T'}^k$, for any sample $\omega\in\Omega$ we get
\begin{equation}\label{nrobust_1}
\xi_{X(0)}(\omega)> kT' \Longleftrightarrow \tau_{T'}^k(\omega)<\infty,
\end{equation}
and if $\tau_{T'}^k(\omega)<\infty$ then
\begin{eqnarray}\label{nrobust_2}
\xi_{X(0)}(\omega)&&=\sum_{s=0}^{\infty} \mathbbm{1}_{\{\mathcal{E}_{X(s+1)}(\omega)\neq \mathcal{E}_{X(s)}(\omega)\}}\nonumber\\
&&=\sum_{s=1}^{\infty} \mathbbm{1}_{\{\mathcal{E}_{X(s)}(\omega)\neq \mathcal{E}_{X(s-1)}(\omega)\}}\nonumber\\
&&=\Big(\sum_{s=1}^{\tau_{T'}^k(\omega)}+\sum_{s=\tau_{T'}^k(\omega)+1}^{\infty}\Big) \mathbbm{1}_{\{\mathcal{E}_{X(s)}(\omega)\neq \mathcal{E}_{X(s-1)}(\omega)\}}\nonumber\\
&&=kT'+1+\xi_{X(\tau_{T'}^k(\omega))}(\omega).
\end{eqnarray}
By (\ref{nrobust_1}), (\ref{nrobust_2}) and Lemma \ref{robust},  for any $k\geq 2$ we get
\begin{eqnarray}\label{nrobust_3}
&&\Prob\left(\xi_{X(0)}> kT'|\xi_{X(0)}>(k-1)T'\right)\nonumber\\
&&~=\Prob\left(\xi_{X(0)}> kT' | \tau_{T'}^{k-1}<\infty \right)\nonumber\\
&&~= \Prob\left(\xi_{X(\tau_{T'}^{k-1})}+(k-1)T'+1> kT'| \tau_{T'}^{k-1}<\infty \right)\nonumber\\
&&~= \Prob\left(\xi_{X(\tau_{T'}^{k-1})}>T'-1| \tau_{T'}^{k-1}<\infty \right)\nonumber\\
&&~\leq  1-\delta,
\end{eqnarray}
where $\delta\in(0,1)$ is the same constant appearing in Lemma \ref{robust}.
By (\ref{nrobust_3}) and Lemma \ref{robust}, for any $k\in\mathbb{Z}^+$ we have
\begin{eqnarray*}
&&\Prob\left(\xi_{X(0)}> kT'\right)\\
&&=\Prob\left(\xi_{X(0)}> kT', \xi_{X(0)}>(k-1)T'\right)\\
&&~~+ \Prob\left(\xi_{X(0)}> kT', \xi_{X(0)}\leq (k-1)T'\right)\\
&&=\Prob\left(\xi_{X(0)}>(k-1)T'\right)\\
&&~~\times\Prob\left(\xi_{X(0)}> kT'|\xi_{X(0)}>(k-1)T'\right)\\
&&\leq (1-\delta) \Prob\left(\xi_{X(0)}>(k-1)T' \right)\\
&&\leq \cdots \leq  (1-\delta)^{k-1} \Prob\left(\xi_{X(0)}> T'\right)\\
&&\leq  (1-\delta)^{k}.
\end{eqnarray*}
Since
\begin{equation*}
\sum_{k=1}^\infty \Prob\left(\xi_{X(0)}> kT'\right)\leq \sum_{k=1}^\infty  (1-\delta)^k<\infty,
\end{equation*}
by the Borel-Cantelli lemma we have
\begin{equation*}
\Prob\left(\bigcap_{m=1}^\infty \bigcup_{k=m}^{\infty}  \left\{\xi_{X(0)}> kT'\right\}\right)
=\Prob\left(\xi_{X(0)}=\infty\right)=0,
\end{equation*}
which indicates  $\Prob\big(\xi_{X(0)}<\infty\big)=1$.

\bibliographystyle{IEEEtran}
\bibliography{alias,Main,FB}

\end{document}